\documentclass[leqno]{article}

\def\diam{\mathop{\rm diam}}
\def\dist{\mathop{\rm dist}}

\begin{document}

\title{Potpourri, 10}

\author{Stephen William Semmes	\\
	Rice University		\\
	Houston, Texas}

\date{}

\maketitle

	Let $(M, d(x, y))$ be a metric space.  Thus $M$ is a nonempty
set and $d(x, y)$ is a real-valued function defined for $x, y \in M$
such that $d(x, y) \ge 0$ for all $x, y \in M$, $d(x, y) = 0$ if and
only if $x = y$, $d(x, y) = d(y, x)$ for all $x, y \in M$, and
\begin{equation}
	d(x, z) \le d(x, y) + d(y, z)
\end{equation}
for all $x, y, z \in M$.  If we have the stronger condition that
\begin{equation}
	d(x, z) \le \max(d(x, y), d(y, z))
\end{equation}
for all $x, y, z \in M$, then we say that $d(x, y)$ is an
\emph{ultrametric} on $M$.

	Of course the real line ${\bf R}$ with the standard metric $|x
- y|$ is a metric space.  We shall discuss some more examples soon.

	Let $(M, d(x, y))$ be a metric space, and let $E$ be a subset
of $M$.  We say that $E$ is bounded if the collection of real numbers
$d(x, y)$, $x, y \in E$ are bounded, and in this event we define the
diameter of $E$ by
\begin{equation}
	\diam E = \sup \{d(x, y) : x, y \in E\}.
\end{equation}
The diameter of the empty set is defined to be $0$, and we can also
define the diameter of an unbounded subset of $M$ to be $+\infty$.

	For each $p \in M$ and $r > 0$ define the open and closed
balls in $M$ with center $p$ and radius $r$ by
\begin{equation}
	B(p, r) = \{x \in M : d(x, p) < r\}
\end{equation}
and
\begin{equation}
	\overline{B}(p, r) = \{x \in M : d(x, p) \le r\},
\end{equation}
respectively.  Thus
\begin{equation}
	B(p, r) \subseteq \overline{B}(p, r)
\end{equation}
and
\begin{equation}
	\diam \overline{B}(p, r) \le 2 \, r
\end{equation}
for all $p \in M$ and $r > 0$.

	If $E$ is a subset of $M$ and $\mathcal{A}$ is a family of
subsets of $M$, then we say that $\mathcal{A}$ is an \emph{admissible
covering} of $E$ if
\begin{equation}
	E \subseteq \bigcup_{A \in \mathcal{A}} A
\end{equation}
and $\mathcal{A}$ has at most countably many elements.  In some
situations one might wish to restrict one's attention to admissible
coverings by bounded subsets of $M$.  One might also be interested in
coverings which contain only finitely many subsets of $M$.

	For example, if $E$ is any subset of $M$, then we can simply
cover $E$ by $E$ itself, or by $M$.  If $p$ is any element of $M$,
then the family of balls $B(p, l)$, with $l$ a positive integer, is an
admissible covering of any subset of $M$ by bounded subsets of $M$.

	Suppose that $h(t)$ is a monotone increasing continuous
real-valued function defined for nonnegative real numbers $t$ such
that $h(0) = 0$, $h(t) > 0$ when $t > 0$.  Put
\begin{equation}
	h(+\infty) = \sup \{h(t) : 0 \le t < \infty\},
\end{equation}
which is equal to $+\infty$ if $h(t)$ is unbounded. 

	If $E$ is any subset of $M$, then we define $\mu_h(E)$, the
Hausdorff content of $E$ relative to the function $h$, to be the
infimum of
\begin{equation}
	\sum_{A \in \mathcal{A}} h(\diam A)
\end{equation}
over all admissible coverings $\mathcal{A}$ of $E$.  As a special
case, if $\alpha$ is a positive real number, then we can take $h(t) =
t^\alpha$.  For this we write $\mu^\alpha(E)$ instead of $\mu_h(E)$,
and $\mu^\alpha(E)$ is called the $\alpha$-dimensional Hausdorff
content of $E$.

	Of course $\mu_h(E) \ge 0$ for all subsets $E$ of $M$, and
$\mu_h(E) = 0$ if $E$ is the empty set, or if $E$ contains only one
element.  Using the trivial covering of $E$ by $E$ we get that
\begin{equation}
\label{mu_h(E) le h(diam E)}
	\mu_h(E) \le h(\diam E).
\end{equation}	

	In general $\mu_h(E)$ may be equal to $+\infty$ when $E$ is
unbounded.  If $h$ is unbounded and we restrict ourselves to
admissible coverings with only finitely many elements, then the
analogue of $\mu_h(E)$ is equal to $+\infty$ when $E$ is unbounded.

	If $E$, $\widetilde{E}$ are subsets of $M$ with
\begin{equation}
	\widetilde{E} \subseteq E,
\end{equation}
then
\begin{equation}
	\mu_h(\widetilde{E}) \le \mu_h(E).
\end{equation}
If $E_1$, $E_2$ are subsets of $M$, then
\begin{equation}
	\mu_h(E_1 \cup E_2) \le \mu_h(E_1) + \mu_h(E_2).
\end{equation}
This follows easily by combining arbitrary admissible coverings of
$E_1$, $E_2$ to get admissible coverings of $e_1 \cup E_2$, and it
would also work for the analogue of $\mu_h(E)$ defined using finite
coverings of $E$.  If $E_1, E_2, E_3, \ldots$ is a sequence of subsets
of $M$, then
\begin{equation}
	\mu_h\bigg(\bigcup_{j=1}^\infty E_j\bigg)
		\le \sum_{j=1}^\infty \mu_h(E_j),
\end{equation}
because arbitrary admissible coverings of the $E_j$'s can be combined
to give an admissible covering of the union.  This uses the fact that
a countable union of countable sets is a countable set.

	For the analogue of $\mu_h$ defined using finite coverings we
have in particular that $\mu_h(E) = 0$ when $E$ is a finite set.  If
we use countable coverings as above then $\mu_h(E) = 0$ when $E$ is at
most countable.

	Notice that the closure of a bounded subset of a metric space
$M$ is also bounded and has the same diameter.  Using this one can
check that one gets the same result for $\mu_h(E)$ if one restricts
one's attention to admissible coverings of $E$ by closed subsets of
$M$.  

	The same remark applies to the analogue of $\mu_h(E)$ based on
finite coverings of $E$.  For this it follows that one gets the same
answer for the closure of $E$ as for $E$.  This does not work in
general for $\mu_h(E)$ based on coverings which are at most countable.

	If $A$ is any subset of $M$ and $r$ is a positive real number,
put
\begin{equation}
	A(r) = \{x \in M : \hbox{ there is an } a \in A
				\hbox{ such that } d(x, a) < r\}.
\end{equation}
Thus $A \subseteq A(r)$, and one can check that $A(r)$ is always an
open subset of $M$.

	If $A$ is a bounded subset of $M$, then $A(r)$ is bounded too,
and one can check that
\begin{equation}
	\diam A(r) \le \diam A + 2 \, r.
\end{equation}
In the case of an ultrametric we have that the diameter of $A(r)$ is
less than or equal to the maximum of the diameter of $A$ and $r$.
Using this one can show that we can restrict ourselves to coverings by
open sets and get the same answers for Hausdorff content.

	If $E$ is a compact subset of $M$, then for each open covering
of $E$ there is a finite subcovering of $E$.  Thus the versions of
$\mu_h(E)$ based on coverings which are finite or at most countable
give the same result when $E$ is compact.

	Let $\epsilon > 0$ be given, and let $E$ be a subset of $M$.
An admissible covering $\mathcal{A}$ of $E$ in $M$ is called an
$\epsilon$-covering if
\begin{equation}
	\diam A < \epsilon
\end{equation}
for all $A \in \mathcal{A}$.  If $M$ is a separable metric space,
meaning that there is a dense subset of $M$ which is at most
countable, then there is an $\epsilon$-covering of $M$, and hence of
any subset of $M$, for all $\epsilon > 0$.

	As a variant of the Hausdorff content associated to $h(t)$,
define $\mathcal{H}_{h, \epsilon}(E)$ to be the infimum of $\sum_{A
\in \mathcal{A}} h(\diam A)$ over all $\epsilon$-coverings
$\mathcal{A}$ of $E$, if there are any, and otherwise put
$\mathcal{H}_{h, \epsilon}(E) = + \infty$.  If $\epsilon_1$,
$\epsilon_2$ are positive real numbers with $\epsilon_1 \le
\epsilon_2$, then
\begin{equation}
	\mathcal{H}_{h, \epsilon_2}(E) \le \mathcal{H}_{h, \epsilon_1}(E).
\end{equation}
For all $\epsilon > 0$ we have that
\begin{equation}
	\mu_h(E) \le \mathcal{H}_{h, \epsilon}(E).
\end{equation}

	Just as for Hausdorff content, we can restrict ourselves to
$\epsilon$-coverings by open or closed subsets of $M$ and get the same
result for $\mathcal{H}_{h, \epsilon}(E)$.  If $E$ is a compact subset
of $M$, then $\mathcal{H}_{h, \epsilon}(E)$ is equal to its analogue
for finite coverings.

	The Hausdorff measure of a subset $E$ of $M$ associated to the
function $h(t)$ is defined by
\begin{equation}
	\mathcal{H}_h(E) 
		= \sup_{\epsilon > 0} \mathcal{H}_{h, \epsilon}(E).
\end{equation}
Thus
\begin{equation}
	\mu_h(E) \le \mathcal{H}_{h, \epsilon}(E) \le \mathcal{H}_h(E)
\end{equation}
for all $\epsilon > 0$.  If $\alpha$ is a positive real number and
$h(t) = t^\alpha$, then we write $\mathcal{H}^\alpha_\epsilon(E)$,
$\mathcal{H}^\alpha(E)$ in place of $\mathcal{H}_{h, \epsilon}(E)$,
$\mathcal{H}_h(E)$.

	Notice that $\mathcal{H}_h(E) = 0$ when $E$ is the empty set
or $E$ consists of just one element.  If $E$, $\widetilde{E}$ are
subsets of $M$ and $\widetilde{E} \subseteq E$, then
$\mathcal{H}_h(\widetilde{E}) \le \mathcal{H}_h(E)$.  

	If $E_1, E_2, \ldots$ are subsets of $M$, then
\begin{equation}
	\mathcal{H}_{h, \epsilon} \bigg(\sum_{j=1}^\infty E_j \bigg)
		\le \sum_{j=1}^\infty \mathcal{H}_{h, \epsilon}(E_j)
\end{equation}
for all $\epsilon > 0$, and therefore
\begin{equation}
	\mathcal{H}_h\bigg(\sum_{j=1}^\infty E_j \bigg)
		\le \sum_{j=1}^\infty \mathcal{H}_h(E_j).
\end{equation}
For the analogue of $\mathcal{H}_h$ based on finite coverings we still
have finite subadditivity.

	Suppose that $E_1$, $E_2$ are subsets of $M$ and that there is
an $\eta > 0$ such that $d(x, y) \ge \eta$ for all $x \in E_1$ and all
$y \in E_2$.  In this event we have that
\begin{equation}
	\mathcal{H}_{h, \epsilon}(E_1 \cup E_2)
	   \ge \mathcal{H}_{h, \epsilon}(E_1) + \mathcal{H}_{h, \epsilon}(E_2)
\end{equation}
for all $\epsilon \in (0, \eta]$.  Hence
\begin{equation}
	\mathcal{H}_h(E_1 \cup E_2)
		\ge \mathcal{H}_h(E_1) + \mathcal{H}_h(E_2).
\end{equation}

	Notice that $\mu_h(E) = 0$ for a subset $E$ of $M$ if and only
if for each $\eta > 0$ there is an admissible covering $\mathcal{A}$
of $E$ such that
\begin{equation}
	\sum_{A \in \mathcal{A}} h(A) < \eta.
\end{equation}
This implies that $\mathcal{H}_h(E) = 0$, since the diameters of the
covering sets $A$ have to be small anyway.

	Let us consider the special case where $M$ is the real line
with the standard metric.  For this we can restrict ourselves to
coverings by intervals in the definitions of Hausdorff measure and
Hausdorff content, since every subset of the real line is contained in
an interval with the same diameter.  We can allow unbounded intervals
such as the real line itself to accommodate unbounded subsets of the
real line.

	If $a$, $b$ are real numbers with $a \le b$, let $[a, b]$ be
the usual closed interval with endpoints $a$, $b$, consisting of the
real numbers $x$ such that $a \le x \le b$.  Recall that this is a
compact subset of the real line.

	Clearly
\begin{equation}
	\mu^1([a, b]) \le b - a.
\end{equation}
For each $\epsilon > 0$ we can cover $[a, b]$ by finitely many
closed intervals, each of length less than $\epsilon$, and so that the
sum of their lengths is equal to $b - a$, which implies that
\begin{equation}
	\mathcal{H}^1_\epsilon ([a, b]) \le b - a
\end{equation}
and therefore
\begin{equation}
	\mathcal{H}^1([a, b]) \le b - a.
\end{equation}

	Conversely,
\begin{equation}
	\mu^1([a, b]) \ge b - a.
\end{equation}
This reduces to the observation that if $I_1, \ldots, I_l$ are
intervals in the real line such that
\begin{equation}
	[a, b] \subseteq I_1 \cup \cdots \cup I_l,
\end{equation}
then $b - a$ is less than or equal to the sum of the lengths of the
$I_j$'s.  In summary we have that
\begin{equation}
	\mathcal{H}^1([a, b]) = \mu^1([a, b]) = b - a.
\end{equation}	

	Let $X_1, X_2, \ldots$ be a sequence of finite sets, where
each $X_j$ has $n_g \ge 2$ elements, and let $X$ denote the space of
sequences $x = \{x_j\}_{j=1}^\infty$ with $x_j \in X_j$ for each $j$.
For each nonnegative integer $l$ and each $x \in X$ let $N_l(x)$ be
the set of $y \in X$ such that the $j$th terms of $x$ and $y$ are
equal when $j \le l$.  These ``neighborhoods'' of points in $X$
generate a topology for $X$ in the usual way, in which a subset $U$ of
$X$ is open if for each $x \in U$ there is an $l \ge 0$ such that
$N_l(x) \subseteq U$.  If $x, y \in X$ and $y \in N_l(x)$, then
$N_l(y) = N_l(x)$, as one can easily see, and thus $N_l(x)$ is an open
subset of $X$.

	Let us call a subset of $X$ of the form $N_l(x)$ for some $x
\in X$ and $l \ge 0$ a cell.  Thus cells are open subsets of $X$, and
one can check that they are also closed, because the complement of a
cell can be expressed as a union of cells.  With this topology $X$
becomes a compact Hausdorff space which is totally disconnected, which
is to say that $X$ does not contain any connected subsets with more
than one element.  If $\mathcal{C}_1$, $\mathcal{C}_2$ are cells in
$X$, then either $\mathcal{C}_1 \subseteq \mathcal{C}_2$, or
$\mathcal{C}_2 \subseteq \mathcal{C}_1$, or $\mathcal{C}_1 \cap
\mathcal{C}_2 = \emptyset$.  We can be more precise and say that if
$x, y \in X$ and $l$, $p$ are nonnegative integers with $l \le p$,
then either $N_p(y) \subseteq N_l(x)$ or $N_l(x) \cap N_p(y) =
\emptyset$.

	Suppose that $\{r_l\}_{l=0}^\infty$ is a strictly decreasing
sequence of positive real numbers such that $\lim_{l \to \infty} r_l =
0$.  Define a distance function on $X$ by saying that the distance
from $x$ to $y$, $x, y \in X$ is equal to $0$ when $x = y$ and to
$r_l$ when the $j$th terms of $x$, $y$ are equal for $j \le l$ and
different for $j = l+1$.  With respect to this distance function,
$N_l(x)$ is the same as the set of $y \in X$ whose distance to $x$ is
less than or equal to $r_l$.  This defines an ultrametric on $X$
compatible with the topology generated by these standard neighborhoods
in $X$.

	If $A$ is a subset of $X$, then there is a cell in $X$ which
contains $A$ and which has the same diameter as $A$, and therefore one
may restrict one's attention to admissible coverings by cells in $X$
in the definition of the Hausdorff content or Hausdorff measure of a
subset of $X$.  Of course
\begin{equation}
	\mu_h(\mathcal{C}) \le h(r_l)
\end{equation}
for any cell $\mathcal{C}$ in $X$ of diameter $r_l$.  If $p > l$, then
$\mathcal{C}$ is the union of
\begin{equation}
	\prod_{i = l + 1}^p n_i
\end{equation}
cells of diameter $r_p$.  Hence
\begin{equation}
	\mathcal{H}_{h, \epsilon}(\mathcal{C}) 
		\le \bigg(\prod_{i=l+1}^p n_i \bigg) \, h(r_p)
\end{equation}
when $r_p < \epsilon$.

	A particularly nice situation occurs when the function $h(t)$
satisfies
\begin{equation}
	h(r_l) = n_l \, h(r_{l+1})
\end{equation}
for all $l \ge 0$.  In this event we get that
\begin{equation}
	\mathcal{H}_h(\mathcal{C}) \le h(r_l)
\end{equation}
for any cell $\mathcal{C}$ in $X$ with diameter $r_l$.  Moreover,
\begin{equation}
	\mu_h(\mathcal{C}) \ge h(r_l)
\end{equation}
under these conditions.  This is because any finite covering of
$\mathcal{C}$ can be refined to a covering of cells of the same
diameter, and any covering of $\mathcal{C}$ by cells of diameter
$r_p$, $p \ge l$, has to include all of the cells contained in
$\mathcal{C}$ of diameter $r_p$.  If the $n_j$'s are all equal to some
positive integer $n$ and $r_l = r^l$ for some $r \in (0, 1)$, then we
can take
\begin{equation}
	h(t) = t^\alpha
\end{equation}
where
\begin{equation}
	n \, r^\alpha = 1.
\end{equation}

	More generally, suppose that $(M, d(x, y))$ is a metric space
with $d(x, y)$ an ultrametric.  Let $p_1$, $p_1$ be elements of $M$,
and let $r_1 \le r_2$ be positive real numbers.  If
\begin{equation}
	B(p_1, r_1) \cap B(p_2, r_2) \ne \emptyset,
\end{equation}
then
\begin{equation}
	B(p_1, r_1) \subseteq B(p_2, r_2).
\end{equation}
Similarly, if
\begin{equation}
	\overline{B}(p_1, r_1) \cap \overline{B}(p_2, r_2) \ne \emptyset,
\end{equation}
then
\begin{equation}
	\overline{B}(p_1, r_1) \subseteq \overline{B}(p_2, r_2).
\end{equation}
The diameter of a ball of radius $r$ in $M$ is less than or equal to
$r$.

	One can check that every open ball in $M$ is a closed subset
of $M$.  Also, every closed ball in $M$ is an open subset of $M$.  In
particular, $M$ is totally disconnected.

	In our space $X$ as before, it is natural to say that a
cell of the form $N_l(x)$ should have measure equal to $1$ when
$l = 0$, and to 
\begin{equation}
	\frac{1}{\prod_{i=1}^l n_i}
\end{equation}
when $l > 0$.  In this way the cells of the same diameter have the
same measure.

	Let $l$ be a nonnegative integer, and let $E_l$ be a subset of
$X$ such that $E_l$ contains exactly one element from each of the
cells in $X$ of diameter $r_l$.  If $f$ is a continuous real-valued
function on $X$, then we can associate to $f$ the Riemann sum
\begin{equation}
	\frac{1}{\prod_{i=1}^l n_i} \sum_{x \in E_l} f(x).
\end{equation}
The integral of $f$ can be defined as a limit of Riemann sums as $l
\to \infty$.  Recall that because $X$ is compact, $f$ is automatically
uniformly continuous on $X$.  This implies that the Riemann sums
converge and are independent of the choices of the $E_l$'s.

	If $\mathcal{C}$ is a cell in $X$, then the function on $X$
equal to $1$ on $\mathcal{C}$ and to $0$ on the complement of
$\mathcal{C}$ in $X$ is a continuous function on $X$.  The integral of
this function is equal to the measure of the cell as defined before.
Basically we simply have a probability distribution on $X$ which is
uniformly distributed at each step in the obvious way.

	If $(M, d(x, y))$ and $(N, \rho(u, v))$ are metric spaces,
then a mapping $f : M \to N$ is said to be $C$-Lipschitz for some $C
\ge 0$ if
\begin{equation}
	\rho(f(x), f(y)) \le C \, d(x, y)
\end{equation}
for all $x, y \in M$.  This is equivalent to saying that for each
bounded subset $A$ of $M$, $f(A)$ is a bounded subset of $N$ with
diameter less than or equal to $C$ times the diameter of $A$ in $M$.
In this case, if $h(t)$ is a function on $[0, \infty)$ as before and
$\widetilde{h}(t) = h(C^{-1} \, t)$, then for each subset $E$ of $M$
we have that the Hausdorff content of $f(E)$ in $N$ with respect to
$\widetilde{h}$ is less than or equal to the Hausdorff content of $E$
in $M$ with respect to $h$.  Similarly, $\mathcal{H}_{\widetilde{h}, C
\, \epsilon}(f(E))$ in $N$ is less than or equal to $\mathcal{H}_{h,
\epsilon}(E)$ in $M$ for all $\epsilon > 0$.  As a consequence,
$\mathcal{H}_{\widetilde{h}}(f(E))$ in $N$ is less than or equal to
$\mathcal{H}_h(E)$ in $M$.

	A real-valued function $f(x)$ on $M$ is $C$-Lipschitz if and
only if
\begin{equation}
	f(x) \le f(y) + C \, d(x, y)
\end{equation}
for all $x, y \in M$.  For each $p \in M$ we have that $f_p(x) = d(x,
p)$ is $1$-Lipschitz.  More generally, if $A$ is a nonempty subset of
$M$, then
\begin{equation}
	\dist(x, A) = \inf \{d(x, y) : y \in A\}
\end{equation}
is $1$-Lipschitz.  If the $1$-dimensional Hausdorff content of $M$ is
equal to $0$, then the same must be true of the image of any
real-valued Lipschitz function on $M$.  In particular, $M$ must be
totally disconnected.

	Let $\phi(t)$ be a monotone increasing continuous real-valued
function defined for nonnegative real numbers $t$ such that $\phi(0) =
0$ and $\phi(t) > 0$ when $t > 0$.  If $(M, d(x, y))$ is a metric
space and $d(x, y)$ is actually an ultrametric on $M$, then $\phi(d(x,
y))$ is also an ultrametric on $M$ which defines the same topology as
$d(x, y)$.  In particular, this works when $\phi(t) = t^a$ for some $a
> 0$.

	Now assume further that $\phi(t)$ is subadditive, so that
\begin{equation}
	\phi(x + y) \le \phi(x) + \phi(y)
\end{equation}
for all $x, y \ge 0$.  This implies that $\phi(t)$ is uniformly
continuous under these conditions.  Namely,
\begin{equation}
	\phi(x) \le \phi(x + r) \le \phi(x) + \phi(r),
\end{equation}
so that uniform continuity follows from continuity at $0$.  If $(M,
d(x, y))$ is a metric space, then $\phi(d(x, y))$ is a metric on $M$
which defines the same topology as $d(x, y)$ does.  If $0 < a \le 1$,
then $\phi(t) = t^a$ satisfies these conditions.

	In either of these two situations, if $A$ is a bounded subset
of $M$ with respect to $d(x, y)$, then $A$ is bounded with respect to
$\phi(d(x, y))$.  The diameter of $A$ with respect to $\phi(d(x, y))$
is equal to $\phi$ of the diameter of $A$ with respect to $d(x, y)$.

\end{document}